\documentclass[smallextended,final]{svjour3}

\usepackage{amssymb,amsfonts,amsmath}
\usepackage[mathscr]{eucal}
\usepackage{stmaryrd}
\usepackage{amsbsy}
\usepackage{bm}
\usepackage{braket}
\usepackage[usenames]{xcolor}
\usepackage{ifpdf}
\usepackage[margin=5pt,justification=centering]{subfig}
\usepackage{array}
\usepackage{environ}
\ifpdf
\usepackage[colorlinks,urlcolor=blue,citecolor=blue,linkcolor=blue]{hyperref}
\else
\newcommand{\href}[2]{{#2}}
\usepackage{url}
\newcommand{\nolinkurl}[1]{\url{#1}}
\fi
\usepackage{xspace}
\usepackage{enumitem}

\usepackage{algorithm}
\usepackage{algpseudocode}
\usepackage{multirow}
\usepackage{mathtools}
\usepackage{tikz}
\usetikzlibrary{calc}

\newcommand{\Tra}{{\sf T}} 
\newcommand{\Parens}[1]{\left(#1\right)}
\newcommand{\Real}{\mathbb{R}}
\newcommand{\qtext}[1]{\quad\text{#1}\quad}
\newcommand{\V}[2][]{{\bm{#1\mathbf{\MakeLowercase{#2}}}}} 
\newcommand{\Vn}[3][]{{\bm{#1\mathbf{\MakeLowercase{#2}}}}^{(#3)}} 
\newcommand{\VnE}[4][]{#1{\MakeLowercase{#2}}^{(#3)}_{#4}} 
\newcommand{\M}[2][]{{\bm{#1\mathbf{\MakeUppercase{#2}}}}} 
\newcommand{\Mn}[3][]{{\bm{#1\mathbf{\MakeUppercase{#2}}}}^{(#3)}} 
\newcommand{\MC}[3][]{\V[#1]{#2}_{#3}} 
\newcommand{\MnC}[4][]{\Vn[#1]{#2}{#3}_{#4}} 
\newcommand{\T}[2][]{\boldsymbol{#1\mathscr{\MakeUppercase{#2}}}} 
\newcommand{\TE}[3][]{#1{\MakeLowercase{#2}}_{#3}}

\definecolor{color1}{rgb}{0.2235,0.4157,0.6941}
\definecolor{color2}{rgb}{0.8549,0.4863,0.1882}
\definecolor{color3}{rgb}{0.2431,0.5882,0.3176}
\definecolor{color4}{rgb}{0.8000,0.1451,0.1608}
\definecolor{color5}{rgb}{0.3255,0.3176,0.3294}
\definecolor{color6}{rgb}{0.4196,0.2980,0.6039}
\definecolor{color7}{rgb}{0.5725,0.1412,0.1569}
\definecolor{color8}{rgb}{0.5804,0.5451,0.2392}
\definecolor{color8}{rgb}{0,0,0}

\ifpdf
\newcommand{\Sec}[1]{\hyperref[sec:#1]{\S\ref*{sec:#1}}} %
\newcommand{\SEC}[1]{\hyperref[sec:#1]{Section~\ref*{sec:#1}}} %
\newcommand{\Eqn}[1]{\hyperref[eq:#1]{{\rm (\ref*{eq:#1})}}} %
\newcommand{\Fig}[1]{\hyperref[fig:#1]{Figure~\ref*{fig:#1}}} %
\newcommand{\Tab}[1]{\hyperref[tab:#1]{Table~\ref*{tab:#1}}} %
\newcommand{\Alg}[1]{\hyperref[alg:#1]{Algorithm~\ref*{alg:#1}}} %
\else
\newcommand{\Sec}[1]{{\S\ref{sec:#1}}} %
\newcommand{\SEC}[1]{{Section~\ref{sec:#1}}} %
\newcommand{\Eqn}[1]{{(\ref{eq:#1})}} %
\newcommand{\Fig}[1]{{Figure~\ref{fig:#1}}} %
\newcommand{\Tab}[1]{{Table~\ref{tab:#1}}} %
\newcommand{\Alg}[1]{{Algorithm~\ref{alg:#1}}} %
\fi

\newcommand{\RTMN}{\ensuremath\mathbb{R}^{[m,n]}}
\newcommand{\STMN}{\ensuremath\mathbb{S}^{[m,n]}}
\newcommand{\TA}{\T{A}}
\newcommand{\MX}{\M{X}}
\newcommand{\Vxk}{\MC{X}{k}}
\newcommand{\Vxkm}{\Vxk^{m}}
\newcommand{\Slk}{\lambda_k}
\newcommand{\Vl}{\V{\lambda}}
\newcommand{\Vi}{\V{i}}
\newcommand{\Vc}{\V{c}}
\newcommand{\FD}[2]{\frac{\partial #1}{\partial #2}}
\newcommand{\Pm}{^{\phantom{m}}}

\begin{document}

\title{Numerical Optimization for Symmetric Tensor Decomposition%
\thanks{This material is based upon work supported by the U.S. Department of Energy, Office of Science, Office of Advanced Scientific Computing Research, Applied Mathematics program. Sandia National Laboratories is a multi-program laboratory managed and operated by Sandia Corporation, a wholly owned subsidiary of Lockheed Martin Corporation, for the U.S. Department of Energy's National Nuclear Security Administration under contract DE-AC04-94AL85000.}}

\author{Tamara G. Kolda}

\institute{Tamara G. Kolda \at 
  Sandia National Laboratories, Livermore, CA\\
  \email{tgkolda@sandia.gov}}

\date{Draft: \today\@ / Received: date / Accepted: date}

\maketitle

\begin{abstract}
We consider the problem of decomposing a real-valued symmetric tensor
as the sum of outer products of real-valued vectors. Algebraic methods
exist for computing complex-valued decompositions of symmetric
tensors, but here we focus on real-valued decompositions, both
unconstrained and nonnegative, for problems with low-rank structure. We discuss when solutions exist and how to formulate the mathematical program.
Numerical results show the properties of the proposed formulations (including one that ignores symmetry) on a set of test problems and illustrate that these straightforward formulations can be effective even though the problem is nonconvex.
\end{abstract}

\section{Introduction}
\label{sec:intro}
We consider the problem of decomposing a real-valued symmetric tensor as the sum of outer products of real-valued vectors.
Let $\TA$ represent an $m$-way, $n$-dimension symmetric tensor. 
Given a real-valued vector $\V{x}$ of length $n$, 
we let $\V{x}^m$ denote the $m$-way, $n$-dimensional 
symmetric outer product tensor such that 
$\Parens{\V{x}^m}_{i_1 i_2 \cdots i_m} = x_{i_1} x_{i_2} \cdots x_{i_m}$.
Comon et~al.~\cite{CoGoLiMo08} showed that
any real-valued symmetric tensor $\TA$  can be decomposed as
\begin{equation}
\label{eq:scp}
  \TA = \sum_{k=1}^p \lambda_k \; \MC{x}{k}^m,
\end{equation}
with $\lambda_k \in \Real$ and $\MC{x}{k} \in \Real^n$ for $k=1,\dots,p$; see the illustration in \Fig{symten}.
We assume that the tensor is low-rank, i.e., $p$ is small relative to
the typical rank of a random tensor.
We survey the methods that have been proposed for related problems and discuss several optimization formulations, including a surprisingly effective method that \emph{ignores} the symmetry. 

\definecolor{TikzCubeColor}{rgb}{0.7,0.7,0.7}
\newcommand{\TikzCube}[3]{ %
  \def\scale{#1};
  \def\xsize{1*\scale};
  \def\ysize{1*\scale};
  \def\xdelta{0.40*\scale};
  \def\ydelta{0.25*\scale};

  \coordinate (FrontLowerLeft) at (#2);
  \coordinate (TopLowerLeft) at ($(FrontLowerLeft)+(0,\ysize)$);
  \coordinate (SideLowerLeft) at ($(FrontLowerLeft)+(\xsize,0)$);
  \coordinate (FrontCenter) at ($(FrontLowerLeft)+(\xsize/2,\ysize/2)$);

  \draw[fill=TikzCubeColor!50] (FrontLowerLeft) -- ++(\xsize,0) -- ++(0,\ysize) -- ++(-\xsize,0) -- cycle;
  \draw[fill=TikzCubeColor!50] (TopLowerLeft)  -- ++(\xdelta,\ydelta) -- ++(\xsize,0) -- ++(-\xdelta,-\ydelta) -- cycle;
  \draw[fill=TikzCubeColor!50] (SideLowerLeft) -- ++(\xdelta,\ydelta) -- ++(0,\ysize) -- ++(-\xdelta,-\ydelta) -- cycle;

  \node at (FrontCenter) {#3};
}
\newcommand{\TikzStar}[3]{
  \def\scale{#1};
  \def\xsize{1*\scale};
  \def\ysize{1*\scale};
  \def\xdelta{0.40*\scale};
  \def\ydelta{0.25*\scale};
  \def\exdelta{0.15*\scale};
  \def\width{0.1*\scale};

  \coordinate (FrontUpperLeft) at ($(#2)+(0,\ysize)$);
  \coordinate (ColUpperLeft) at ($(FrontUpperLeft)+(0,-\exdelta/4)$);
  \coordinate (TubeLowerLeft) at ($(FrontUpperLeft)+(0,\exdelta/4)$);
  \coordinate (RowUpperLeft) at ($(FrontUpperLeft)+(\exdelta,0)$);

  \draw[fill=TikzCubeColor!50] (RowUpperLeft) node[left=2] {$\lambda_{#3}$} -- ++(0,-\width) -- ++(\xsize,0) -- ++(0,\width) node[right] {$\MC{x}{#3}$} -- cycle;
  \draw[fill=TikzCubeColor!50] (TubeLowerLeft) -- ++(1.25*\width,0) -- ++(\xdelta,\ydelta) -- ++(-1.25*\width,0) node[above] {$\MC{x}{#3}$} -- cycle;
  \draw[fill=TikzCubeColor!50] (ColUpperLeft) -- ++(\width,0) -- ++(0,-\ysize) -- ++(-\width,0) node[right=2] {$\MC{x}{#3}$} -- cycle;
 
}

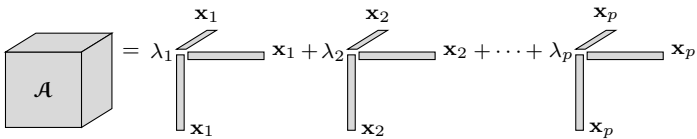
\begin{figure}[htbp]
  \centering
  \begin{tikzpicture}
    \TikzCube{1}{0,0}{$\T{A}$};
    \node at (1.65,1) {$=$};
    \TikzStar{1}{2.25,0}{1};
    \node at (4.0,1) {$+$};
    \TikzStar{1}{4.5,0}{2};
    \node[right] at (6.05,1) {$+\cdots+$};
    \TikzStar{1}{7.5,0}{p};   
  \end{tikzpicture}
  \caption{Symmetric tensor factorization for $m=3$.}
  \label{fig:symten}
\end{figure}

We also consider the related problem of decomposing a real-valued
nonnegative symmetric tensor as the sum of outer products of
real-valued nonnegative vectors.
Let $\TA \geq 0$ represent an $m$-way, $n$-dimension nonnegative symmetric tensor. 
In this case, the goal is a factorization of the form
\begin{equation}
\label{eq:nnscp}
  \TA = \sum_{k=1}^p \MC{x}{k}^m 
  \qtext{with} \MC{x}{k} \geq 0.
\end{equation}
If such a factorization exists, 
we say that $\TA$ is completely positive \cite{QiXuXu13}.
If such a factorization does not exist, then we propose to solve a
``best fit'' problem instead.

The paper is structured as follows.
\SEC{background} provides notation and background material.
Related decompositions, including the best symmetric rank-1 approximation, the symmetric Tucker decomposition, partially symmetric decompositions, and the complex-valued canonical decompositions are discussed in \SEC{related}.
We describe two optimization formulations for symmetric decomposition in \SEC{numerical}, and a mathematical program for the nonnegative problem in \SEC{nonnegative}.
Numerical results, including the methodology for generating challenging problems, is presented in \SEC{results}. 
Finally, \SEC{conclusions} discusses our findings and future challenges.

\section{Background}
\label{sec:background}

\subsection{Notation and preliminaries}

A tensor is a multidimensional array. The number of ways or modes is
called the \emph{order} of a tensor. For example, a matrix is a tensor
of order two. Tensors of order three or greater are called
\emph{higher-order} tensors.

Let $n_1 \times n_2 \times \cdots \times n_m$ denote the size of an
$m$-way tensor. We say that the tensor is \emph{cubic} if all the
modes have the same size, i.e., $n = n_1 = n_2 \cdots = n_m$.
In other words, ``cubic'' is the tensor generalization of ``square.''
In this case, we refer to $n$ as the \emph{dimension} of the tensor.
We let $\Real^{[m,n]}$ denote the space of all
cubic real-valued tensors of order $m$ and dimension $n$.
As appropriate, we use \emph{multiindex} notation to compactly index
tensors so that $\Vi = (i_1,i_2,\dots,i_m)$. Thus,
$a_{\Vi}$ denotes $a_{i_1 i_2 \cdots i_m}$.

The norm of a tensor $\T{A} \in \Real^{[m,n]}$ is the square root of
the sum of squares of its elements, i.e.,
\begin{displaymath}
  \| \TA \| = \sqrt{\sum_{i_1=1}^n \sum_{i_2=1}^n \cdots \sum_{i_m=1}^n a_{\Vi}^2}.
\end{displaymath}
Unless otherwise noted, all norms are the (elementwise) $\ell_2$-norm.

\subsection{Symmetric tensors}

A tensor is symmetric if its entries do not change under permutation of the indices.
Formally, 
we let $\pi(m)$ denote the set of permutations of length $m$. For instance, 
\begin{displaymath}
  \pi(3) = \Set{(1,2,3), (1,3,2), (2,1,3), (2,3,1), (3,1,2), (3,2,1)}.
\end{displaymath}
It is well known that $|\pi(m)| = m!$.
We say a real-valued $m$-way $n$-dimensional tensor 
$\TA$ is \emph{symmetric} \cite{CoGoLiMo08} if
  \begin{displaymath}
    \TE{a}{i_{p(1)} \cdots i_{p(m)}} = \TE{a}{i_1 \cdots i_m}
    \qtext{for all} i_1, \dots, i_m \in \set{1,\dots,n}
    \text{ and } p \in \pi(m).
  \end{displaymath}
Such tensors are also sometimes referred to as \emph{supersymmetric}.
For a 3-way tensor $\T{A}$ of dimension $n$, symmetry means
\begin{displaymath}
  a_{ijk} = a_{ikj} = a_{jik} = a_{kij} = a_{jki} = a_{kji}
  \qtext{for all} i,j,k \in \set{1,\dots,n}.
\end{displaymath}
We let $\STMN \subset \RTMN$ denote
the subspace of all symmetric tensors.

\subsection{Symmetric outer product tensors}

A tensor in $\STMN$ is called \emph{rank one} if it has the form
$\lambda \V{x}^m$ where $\lambda \in \Real$  and $\V{x} \in
\Real^n$.
If $m$ is odd or $\lambda > 0$, then the $m$th real root of $\lambda$ always
exists, so we can rewrite the tensor as
\begin{equation}
  \label{eq:absorb}
  \lambda \V{x}^m = \V{y}^m \qtext{where} \V{y} = \Parens{\sqrt[m]{\lambda}} \V{x}.
\end{equation}
If $m$ is even, however, the $m$th real root does not exist if
$\lambda < 0$, so the scalar cannot be absorbed as in \Eqn{absorb}.

\subsection{Model parameters}
For the symmetric decomposition, we let $\V{\lambda}$ denote the vector of weights and $\M{X}$ denote the matrix of component vectors, i.e.,
\begin{displaymath}
  \V{\lambda} = 
  \begin{bmatrix}
    \lambda_1 & \lambda_2 & \cdots \lambda_p
  \end{bmatrix}^{\Tra}
  \qtext{and}
  \M{X} = 
  \begin{bmatrix}
    \MC{x}{1} & \MC{x}{2} & \cdots & \MC{x}{p}
  \end{bmatrix}.
\end{displaymath}
The notation $x_{ik}$ refers to the $i$th entry in the $k$th column,
so recalling the multiindex notation $\V{i} = (i_1,\dots,i_m)$, we have
\begin{displaymath}
  (\Vxkm)_{\Vi} = x_{i_1 k} x_{i_2 k} \cdots x_{i_m k}.
\end{displaymath}

\section{Related problems}
\label{sec:related}

\subsection{Canonical polyadic  tensor decomposition}

Canonical polyadic (CP) tensor decomposition has been known since 1927 \cite{Hi27,Hi27a}. It is known under several names, two of the most prominent being CANDECOMP as proposed by Carroll and Chang \cite{CaCh70} and PARAFAC by Harshman \cite{Ha70}. Originally, the term CP was proposed as a combination of these two names \cite{Ki00}, but more recently has been re-purposed to mean ``canonical polyadic.'' 
For details on CP, we refer the reader to the survey \cite{KoBa09}. Here, we describe the problem in the case of a cubic tensor $\TA \in \Real^{[m,n]}$. Our goal is to discover a decomposition of the form
\begin{equation}\label{eq:cp}
  \T{A} = \sum_{k=1}^p \MnC{u}{1}{k} \circ \MnC{u}{2}{k} \circ \cdots
  \circ \MnC{u}{m}{k}.
\end{equation}
The circle denotes the vector outer product so the $\Vi = (i_1,i_2,\dots,i_m)$ entry is
\begin{displaymath}
  \Parens{ \Vn{u}{1} \circ \Vn{u}{2} \circ \cdots \circ \Vn{u}{m} }_{\Vi} = 
  \VnE{u}{1}{i_1} \VnE{u}{2}{i_2} \cdots \VnE{u}{m}{i_m}.
\end{displaymath}
Each summand is called a \emph{component}. An illustration is shown in \Fig{cp}.
One of the most effective methods for this problem is alternating least squares. We solve for each \emph{factor} matrix
\begin{displaymath}
  \Mn{U}{j} = 
  \begin{bmatrix}
    \MnC{u}{j}{1} & \MnC{u}{j}{2} & \cdots & \MnC{u}{j}{p}
  \end{bmatrix},
\end{displaymath}
in turn by solving a linear least squares problem, cycling through
all modes (i.e., $j=1,\dots,m$) repeatedly until convergence. See, e.g., \cite[Figure 3.3]{KoBa09} for details.

\newcommand{\TikzStarNS}[3]{
  \def\scale{#1};
  \def\xsize{1*\scale};
  \def\ysize{1*\scale};
  \def\xdelta{0.40*\scale};
  \def\ydelta{0.25*\scale};
  \def\exdelta{0.15*\scale};
  \def\width{0.1*\scale};

  \coordinate (FrontUpperLeft) at ($(#2)+(0,\ysize)$);
  \coordinate (ColUpperLeft) at ($(FrontUpperLeft)+(0,-\exdelta/4)$);
  \coordinate (TubeLowerLeft) at ($(FrontUpperLeft)+(0,\exdelta/4)$);
  \coordinate (RowUpperLeft) at ($(FrontUpperLeft)+(\exdelta,0)$);

  \draw[fill=TikzCubeColor!50] (RowUpperLeft)  -- ++(0,-\width) -- ++(\xsize,0) -- ++(0,\width) node[right] {$\MnC{u}{2}{#3}$} -- cycle;
  \draw[fill=TikzCubeColor!50] (TubeLowerLeft) -- ++(1.25*\width,0) -- ++(\xdelta,\ydelta) -- ++(-1.25*\width,0) node[above] {$\MnC{u}{3}{#3}$} -- cycle;
  \draw[fill=TikzCubeColor!50] (ColUpperLeft) -- ++(\width,0) -- ++(0,-\ysize) -- ++(-\width,0) node[right=2] {$\MnC{u}{1}{#3}$} -- cycle;
 
}

\begin{figure}[htbp]
  \centering
  \begin{tikzpicture}
    \TikzCube{1}{0,0}{$\T{A}$};
    \node at (1.65,1) {$=$};
    \TikzStarNS{1}{2.00,0}{1};
    \node at (4.0,1) {$+$};
    \TikzStarNS{1}{4.25,0}{2};
    \node[right] at (6.05,1) {$+\cdots+$};
    \TikzStarNS{1}{7.25,0}{p};   
  \end{tikzpicture}
  \caption{CP tensor factorization for $m=3$.}
  \label{fig:cp}
\end{figure}
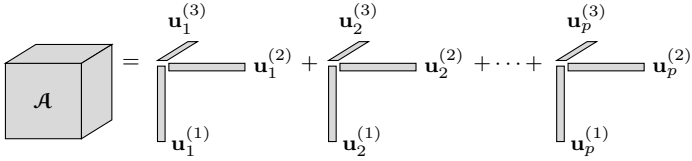

\subsection{Canonical decomposition with partial symmetry}

Partial symmetry has been considered since the work of Carroll and Chang~\cite{CaCh70}. At the same time Carroll and Chang~\cite{CaCh70} introduced CANDECOMP, they also defined INDSCAL which assumes two modes are symmetric. For simplicity of discussion, we assume a cubic tensor $\T{A} \in \Real^{[m,n]}$. For $m=3$ and the last two dimensions being symmetric, this means
\begin{displaymath}
  a_{ijk} = a_{ikj} \qtext{for all} i,j,k \in \set{1,\dots,n},
\end{displaymath}
and the factorization should be of the form
\begin{displaymath}
  \T{A} = \sum_{k=1}^p \MC{u}{k} 
  \circ \MC{v}{k} \circ \MC{v}{k}.
\end{displaymath}
In other words, the last two vectors in each component are equal.
An illustration is provided in \Fig{indscal}.
\newcommand{\TikzStarI}[3]{
  \def\scale{#1};
  \def\xsize{1*\scale};
  \def\ysize{1*\scale};
  \def\xdelta{0.40*\scale};
  \def\ydelta{0.25*\scale};
  \def\exdelta{0.15*\scale};
  \def\width{0.1*\scale};

  \coordinate (FrontUpperLeft) at ($(#2)+(0,\ysize)$);
  \coordinate (ColUpperLeft) at ($(FrontUpperLeft)+(0,-\exdelta/4)$);
  \coordinate (TubeLowerLeft) at ($(FrontUpperLeft)+(0,\exdelta/4)$);
  \coordinate (RowUpperLeft) at ($(FrontUpperLeft)+(\exdelta,0)$);

  \draw[fill=TikzCubeColor!50] (RowUpperLeft)  -- ++(0,-\width) -- ++(\xsize,0) -- ++(0,\width) node[right] {$\MC{v}{#3}$} -- cycle;
  \draw[fill=TikzCubeColor!50] (TubeLowerLeft) -- ++(1.25*\width,0) -- ++(\xdelta,\ydelta) -- ++(-1.25*\width,0) node[above] {$\MC{v}{#3}$} -- cycle;
  \draw[fill=TikzCubeColor!50] (ColUpperLeft) -- ++(\width,0) -- ++(0,-\ysize) -- ++(-\width,0) node[right=2] {$\MC{u}{#3}$} -- cycle;
 
}

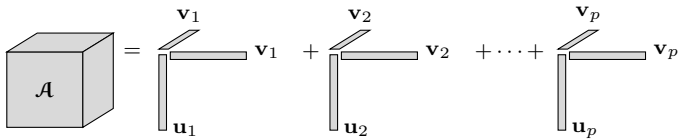
\begin{figure}[htbp]
  \centering
  \begin{tikzpicture}
    \TikzCube{1}{0,0}{$\T{A}$};
    \node at (1.65,1) {$=$};
    \TikzStarI{1}{2.00,0}{1};
    \node at (4.0,1) {$+$};
    \TikzStarI{1}{4.25,0}{2};
    \node[right] at (6.05,1) {$+\cdots+$};
    \TikzStarI{1}{7.25,0}{p};   
  \end{tikzpicture}
  \caption{INDSCAL tensor factorization for $m=3$.}
  \label{fig:indscal}
\end{figure}

Carroll and Chang~\cite{CaCh70} proposed to use an alternating method that ignores symmetry, with the idea that it will often converge to a symmetric solution (up to diagonal scaling).
Later work showed that not all KKT points satisfy this condition \cite{DoTe08}. In \Sec{nosym}, we show how a generalization of this method can be surprisingly effective for symmetric tensor decomposition and provide some motivation for why this might be the case.

We also note that the methods proposed in this manuscript can be extended to partial symmetries.

\subsection{Best symmetric rank-1 approximation}

\begin{figure}[htbp]
  \centering
  \begin{tikzpicture}
    \TikzCube{1}{0,0}{$\T{A}$};
    \node at (1.8,1) {$\approx$};
    \TikzStar{1}{2.30,0}{};   
  \end{tikzpicture}
  \caption{Best symmetric rank-1 decomposition for $m=3$.}
  \label{fig:rank1}
\end{figure}
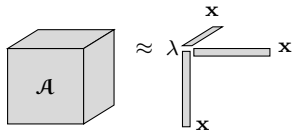

The best symmetric rank-1 approximation problem is
\begin{equation}
  \label{eq:br1}
  \min \left\| \T{A} - \lambda \V{x}^m \right\|^2
  \qtext{subject to} 
  \lambda \in \Real, \; \V{x} \in \Real^n.
\end{equation}
An illustration is shown in \Fig{rank1}.
This problem was first considered in De Lathauwer et al.~\cite{DeDeVa00a}, but their proposed symmetric higher-order power method was not convergent. The power method has been improved so that it is convergent in subsequent work
\cite{KoRe02,ReKo03,KoMa11,GTEP-arXiv-1401.1183}.

This problem is directly related to the problem of computing tensor Z-eigenpairs.  A pair $(\lambda, \V{x})$ is a Z-eigenpair \cite{Li05,Qi05} of a tensor $\TA \in \STMN$ if 
\begin{displaymath}
  \TA \V{x}^{m-1} = \lambda \V{x} \qtext{and} \|\V{x}\|=1,
\end{displaymath}
where $\TA \V{x}^{m-1}$ denotes a vector in $\Real^n$ such that 
\begin{displaymath}
  \Parens{ \TA \V{x}^{m-1} }_{i_1} = \sum_{i_2} \cdots \sum_{i_m} a_{\Vi} x_{i_2} \cdots x_{i_m}
  \qtext{for}
  i_1 \in \set{1,\dots,n}.
\end{displaymath}
The problems are related because any Karush-Kuhn-Tucker (KKT) point of \Eqn{br1} is a Z-eigenpair of $\TA$; see, e.g., \cite{KoMa11}.

Han \cite{Ha12} has considered an unconstrained optimization
formulation of the problem \Eqn{br1}. Cui, Dai, and Nie
\cite{CuDaNi14} use Jacobian SDP relaxations in polynomial
optimization to find \emph{all} real eigenvalues sequentially, from
the largest to the smallest. Nie and Wang \cite{NiWa14} consider
semidefinite relaxations.

\subsection{Symmetric Tucker decomposition}

\begin{figure}[htbp]
  \centering
  \begin{tikzpicture}
    \TikzCube{1}{0,0}{$\T{A}$};
    \node at (1.8,1) {$\approx$};
    \draw[fill=TikzCubeColor!50] (2.25,0) rectangle ++ (0.5,1); 
    \node at (2.5,0.5) {$\M{U}$};
    \draw[fill=TikzCubeColor!50] (4,0.5) rectangle ++ (1,0.5);
    \node at (4.5,0.75) {$\M{U}$};
    \draw[fill=TikzCubeColor!50] (3.25,1.25) -- ++(.5,.35) -- ++(0.5,0) -- ++(-.5,-.35) -- cycle;
    \TikzCube{0.5}{3,0.5}{$\T{B}$};
    \node at (3.7,1.4) {$\M{U}$};
  \end{tikzpicture}
  \caption{Symmetric Tucker decomposition for $m=3$.}
  \label{fig:symtucker}
\end{figure}
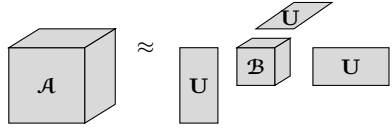

A related problem is symmetric Tucker decomposition. Here the goal is to find an orthogonal matrix $\M{U} \in \Real^{n \times p}$ and a symmetric tensor $\T{B} \in \mathbb{S}^{[m,p]}$ that solves
\begin{displaymath}
  \min \left\| \T{A} - \T[\hat]{A} \right \|^2
  \text{ subject to }
  \TE[\hat]{a}{\Vi} = \sum_{j_1=1}^p \sum_{j_2=1}^p \cdots \sum_{j_m=1}^p
  \TE{b}{j_1j_2\cdots j_m} \; u_{i_1 j_1} u_{i_2 j_2} \cdots u_{i_m j_m}
  .
\end{displaymath}
An illustration is shown in \Fig{symtucker}.
This topic has been considered in
\cite{CaDeDe99,Re13,IsAbVa13} and is useful for compression and signal
processing applications. Alas,  the computational techniques are quite
different, so we do not consider them further in this work.

\subsection{Complex-valued symmetric tensor decomposition}
\label{sec:theoretical}
An alternative version of the problem allows a complex decomposition, i.e., 
\begin{equation}
\label{eq:cscp}
  \TA = \sum_{k=1}^p \MC{x}{k}^m 
  \qtext{with} \MC{x}{k} \in \mathbb{C}^n \text{ for } k = 1,\dots,p.
\end{equation}
Techniques from algebraic geometry have been proposed to solve
\Eqn{cscp} in \cite{BrCoMoTs10,BaBe11,BeGiId11,OeOt13}. More recently, Nie~\cite{Ni14}
devised has a combination of algebraic and numerical approaches for
solving this problem. Generally, these approaches do not scale to
large $n$, though Nie's numerical method scales much better than previous approaches.

In the complex case, the typical rank (i.e., with probability one) is given by the theorem below.
To the best of our knowledge, for the real case, no analogous results are known \cite{CoGoLiMo08}. 

\begin{theorem}[Alexander-Hirschowitz \cite{AlHi95,CoGoLiMo08}]
  For $m > 2$, the typical symmetric rank (over $\mathbb{C}$) of an order-$m$ symmetric tensor of dimension $n$ is
  \begin{displaymath}
    \left\lceil
      \frac{1}{n} \binom{n+k-1}{k}
    \right\rceil
  \end{displaymath}
except for $(m,n) \in \set{(3,5),(4,3),(4,4),(4,5)}$ where it should be increased by one.
\end{theorem}

\section{Optimization formulations for symmetric tensor decomposition}
\label{sec:numerical}

\subsection{Index multiplicities}

A tensor $\TA\in \STMN$ has $n^m$ entries, but not all are distinct. 
Let the set of all possible indices be denoted by
\begin{displaymath}
  \mathcal{R} = \set{ (i_1,\dots, i_m) | i_1,\dots,i_m \in \set{1,\dots,n}}.
\end{displaymath}
Clearly, $|\mathcal{R}| = n^m$.

Following \cite{BaKoPl11}, we define an \emph{index class} as a set of tensor indices such that the corresponding tensor entries all share a value due to symmetry.  For example, for $m=3$ and $n=2$, the tensor indices $(1,1,2)$ and $(1,2,1)$ are in the same index class since $a_{112}=a_{121}$. 
For each index class, we specify an \emph{index representation} which is an index such that the entries are in nondecreasing order.  For instance, $(1,1,2)$ is the index representation for the index class that includes $a_{121}$.
The set
\begin{displaymath}
  \mathcal{I} = \Set{ (i_1,\dots,i_m) | i_1, \dots, i_m \in  \set{1,\dots,n} \text{ and } i_1 \leq i_2 \leq \dots \leq i_m} \subset \mathcal{R}
\end{displaymath}
denotes all possible index representations. 

Each index class also has a \emph{monomial representation} \cite{BaKoPl11}.
For each $\Vi \in \mathcal{I}$ there is a corresponding monomial representation $\Vc$ such that 
\begin{displaymath}
  x_{i_1} x_{i_2} \cdots x_{i_m},
  = x_1^{c_1} x_2^{c_2} \cdots x_n^{c_n}.
\end{displaymath}
Specifically, $c_j$ represents that number of occurrences of index $j$ in $\Vi$ for $j=1,\dots,n$. Clearly, $\sum_j c_j = m$.
Conversely, for a given $\Vc$, we build an index $\Vi$ with $c_1$ copies of 1, $c_2$ copies of 2, etc. This results in an $m$-long index representation.
The set of monomial representations is denoted by
\begin{displaymath}
  \mathcal{C} = \Set{ (c_1,\dots,c_n) | c_1, \dots, c_n \in  \set{0,\dots,m} \text{ and } c_1 + \cdots + c_n = m }.
\end{displaymath}

From \cite{BaKoPl11},
we have that the number of distinct entries of $\T{A}$ is given by
\begin{displaymath}
  |\mathcal{I}| = |\mathcal{C}| = \binom{m+n-1}{m} = \frac{n^m}{m!} + O(n^{m-1}).
\end{displaymath}

It can be shown \cite{BaKoPl11} that the multiplicity of the entry corresponding to a monomial representation $\Vc$ is
\begin{equation}\label{eq:multiplicity}
  \sigma_{\Vc} = 
  \binom{m}{c_1, c_2, \cdots, c_n} = 
  \frac{m!}{c_1! \; c_2! \cdots c_n!} .
\end{equation}
\Tab{monomial} shows an example of index and monomial representations for $\mathbb{S}^{[3,2]}$, including the multiplicities of each element.

\begin{table}[thbp]
  \centering
  \begin{tabular}{|c|c|c|}
    \hline
    Index ($\mathcal{I}$) & Monomimal ($\mathcal{C}$) & Multiplicity ($\sigma$)\\  \hline
    (1,1,1) & (3,0) & 1 \\ \hline
    (1,1,2) & (2,1) & 3 \\ \hline
    (1,2,2) & (1,2) & 3 \\ \hline
    (2,2,2) & (0,3) & 1 \\ \hline
 \end{tabular}
  \caption{Index and monomial representations for $\mathbb{S}^{[3,2]}$.}
  \label{tab:monomial}
\end{table}

Without loss of generality, we exploit the one-to-one correspondence between index and monomial representations to change between them. For example,
\begin{displaymath}
  \| \TA \|^2
  = \sum_{\Vi \in \mathcal{R}} a_{\Vi}^2
  = \sum_{\Vi \in \mathcal{I}} \sigma_{\Vi} a_{\Vi}^2
  = \sum_{\Vc \in \mathcal{C}} \sigma_{\Vc} a_{\Vc}^2,
\end{displaymath}
and
\begin{displaymath}
  (\Vxkm)_{\Vi} = x_{i_1 k} x_{i_2 k} \cdots x_{i_m k} = 
  (\Vxkm)_{\Vc} = x_{1 k}^{c_1} x_{2 k}^{c_2} \cdots x_{n k}^{c_n}. 
\end{displaymath}

\subsection{Two formulations}

For given $\TA \in \STMN$ and $p$, our goal is to find $\Vl$ and $\MX$ such
that \Eqn{scp} is satisfied in a minimization sense. 
We consider two optimization formulations.
The first formulation is the standard least squares formulation, i.e.,
\begin{equation}\label{eq:f1}
  f_1 (\Vl, \MX) 
  = \sum_{\Vi \in \mathcal{R}}
  \left( a_{\Vi} - \sum_{k=1}^p \Slk\Pm (\Vxkm)_{\Vi} \right)^2
  = \sum_{\Vi \in \mathcal{I}}
  \sigma_{\Vi}
  \left( a_{\Vi} - \sum_{k=1}^p \Slk\Pm (\Vxkm)_{\Vi} \right)^2.
\end{equation}
Observe that this counts each unique entry multiple times, according to its multiplicity. The second formulation counts each unique entry only once, i.e.,
\begin{equation}\label{eq:f2}
  f_2 (\Vl, \MX) 
  = \sum_{\Vi \in \mathcal{I}}
  \left( a_{\Vi} - \sum_{k=1}^p \Slk (\Vxkm)_{\Vi} \right)^2.
\end{equation}
Either formulation can be expressed generically as
\begin{displaymath}
  f_{\V{w}}(\Vl,\MX) = \sum_{\Vi \in \mathcal{I}} w_{\Vi} \left( a_{\Vi} - \sum_{k=1}^p \Slk (\Vxkm)_{\Vi} \right)^2
  = \sum_{\Vi \in \mathcal{I}} w_{\Vi} \delta_{\Vi}^2.
\end{displaymath}
Choosing $w_{\Vi} = \sigma_{\Vi}$ yields $f_1$ whereas $w_{\Vi} = 1$ yields $f_2$.
The value $\delta_{\Vi}$ denotes the difference between the model and the tensor at entry $\Vi$.
Note that this formulation easily adapts to the case of missing data, i.e., missing data should have weight of zero in the optimization formulation \cite{AcDuKoMo10,AcDuKoMo11}.

\subsection{Gradients}
\label{sec:gradients}

Using the generic formulation, the gradients are given by
\begin{equation}
  \label{eq:grad}
\begin{aligned}
  \FD{f_{\V{w}}}{\Slk} & = -2\sum_{\Vi \in \mathcal{I}}
  w_{\Vi} \, \delta_{\Vi} \,
  (\Vxkm)_{\Vi}, \\
  \FD{f_{\V{w}}}{x_{jk}} & = -2 \lambda_k \sum_{\Vc \in \mathcal{C}} 
  c_j \, w_{\Vc} \, \delta_{\Vc} \,
  \left( x_{1k}^{c_1} \cdots x_{jk}^{c_j-1} \cdots x_{nk}^{c_n} \right).
\end{aligned}
\end{equation}

For $f_1$, we mention an alternate gradient expression because it is
more efficient to compute for larger values $n$ and $m$. 
The derivation follows \cite{AcDuKo11}, and the gradients are given by
\begin{equation}\label{eq:f1-grad-alt}
\begin{aligned}
  \FD{f_1}{\Slk} & = -2 \TA\Vxkm + 2 \sum_{\ell=1}^p \lambda_{\ell} \left( \Vxk^{\Tra} \V{x}_{\ell} \right)^m, \\
  \FD{f_1}{\Vxk} & = -2m \Slk \TA \Vxk^{m-1} + 2 m \Slk \sum_{\ell=1}^p \lambda_{\ell} \left( \Vxk^{\Tra} \V{x}_{\ell} \right)^{m-1} \V{x}_{\ell}.
\end{aligned}
\end{equation}
This formulation does not easily accommodate missing data since $\V{w}$
is implicit.
\subsection{Scaling ambiguity}
Observe that either objective function suffers from scaling ambiguity.
Suppose we have two equivalent models defined by
\begin{displaymath}
  \sum_{k=1}^p \Slk\Pm \Vxkm
  =
  \sum_{k=1}^p \hat\lambda_k\Pm \V[\hat]{x}_k^m,
\end{displaymath}
related by a positive scaling vector $\V{\rho} \in \Real_+^p$ such that
\begin{displaymath}
  \hat\lambda_k = \rho_k^m \lambda_k\Pm
  \qtext{and}
  \V[\hat]{x}_k =  \V{x}_k / \rho_k
  \qtext{for} k = 1,\dots,p.
\end{displaymath}
To avoid this ambiguity, it is convenient to require $\|\Vxk\| = 1$
for all $k$. We could enforce this condition as an equality
constraint, but instead we treat it as a exact penalty, i.e.,
\begin{equation}
  \label{eq:gamma}
  p_\gamma(\MX) = \gamma \sum_{k=1}^p \Parens{ \Vxk^{\Tra}\Vxk - 1 }^2.
\end{equation}
It is straightforward to observe that the gradient is given by
\begin{displaymath}
  \FD{p_{\gamma}}{\Vxk} = 4 \gamma \Parens{ \Vxk^{\Tra}\Vxk - 1 } \Vxk.
\end{displaymath}
In the experimental results, we see that choosing $\gamma=0.1$ appears
to be adequate for enforcing the penalty.

\subsection{Sparse component weights}

We assume so far that $p$ is known, but this is not always the
case. One technique to get around this problem is to guess a large
value for $p$ and then add a sparsity penalty on $\Vl$, the weight
vector. Specifically, we can use an approximate $\ell_1$ penalty of
the form suggested by \cite{ScFuRo07}:
\begin{displaymath}
  p_{\alpha,\beta} (\Vl)
  = \frac{\beta}{\alpha}
  \sum_{k=1}^p
  \log(1+\exp(-\alpha \lambda_k)) + \log(1 + \exp(\alpha \lambda_k))
  \approx \beta \| \Vl \|_1
\end{displaymath}
In this case, the gradient is
\begin{displaymath}
  \FD{p_{\alpha,\beta} }{\lambda_k}
  = {\beta}
  \left[ 
   (1+\exp(-\alpha \lambda_k))^{-1} + (1 + \exp(\alpha \lambda_k))^{-1}
   \right]
   .
\end{displaymath}
Note that the $\beta$ term is not part of the approximation but rather
the weight of the penalization. In our experiments, the results are
insensitive to the precise choices of $\alpha$ and $\beta$.
\subsection{Putting it all together}
The final function to be optimized is
\begin{displaymath}
  \hat f(\Vl,\MX) = f_{\V{w}}(\Vl,\MX) + p_{\gamma}(\MX) + p_{\alpha,\beta}(\Vl).
\end{displaymath}
The choice of $\V{w}$ determines the choice of objective function. We can also set $w_{\Vi}=0$ for any missing values. The choice of $\gamma$ determines the weight of the penalty on the norm of the columns of $\MX$. 
Since this constraint is easy to satisfy and mostly convenience, the
exact choice of $\gamma$ is not critical. We later show experiments
with $\gamma =0$ and $\gamma = 0.1$, to contrast the difference
between no penalty and a small penalty. (Increasing $\gamma$ beyond
$0.1$ did not have any impact on the experiments.)
The parameter $\alpha$ determines the ``steepness'' of the approximate
$\ell_1$ penalty function, and the choice of $\beta$ determines the
weight of the sparsity-encouraging penalty. In \cite{ScFuRo07}, they
start with a small value of $\alpha$ and gradually increase it.  
In our experiments, we use fixed values $\alpha=10$. The $\beta$ term
is the weight given to the penalty, which is usually determined
heuristically; we use $\beta=0.1$ in our experiments. 

\subsection{Ignoring symmetry}
\label{sec:nosym}

Another approach to symmetric decomposition is to ignore the symmetry altogether and use a standard CP tensor decomposition method such as alternating least squares (ALS) \cite{FaBrHo03,KoBa09}; surprisingly, there are situations under which this non-symmetric method yields a symmetric solution.

Under mild conditions, the CP
decomposition \Eqn{cp} is unique up to permutation and scaling of the components, i.e., \emph{essentially unique}. 
Sidiropoulos and Bro \cite[Theorem 3]{SiBr00} have a general \emph{a posteriori} result on the
essential uniqueness of the CP decompositions for tensors. If we
specialize this result to the symmetric case by assuming $\Mn{U}{j} = \MX$ for $j=1,\dots,m$, the result says that a
sufficient condition for the uniqueness of \Eqn{cp} is
\begin{equation}
  \label{eq:unique}
  2p + (m-1) \leq m \, \text{k-rank}(\M{X}).
\end{equation}
Here, the k-rank of the matrix $\M{X}$ is the largest number $k$
such that every subset of $k$ columns of $\M{X}$ is linearly independent.
\Tab{krank} shows sufficient k-rank's for various values of $m$ and
$p$. For instance, if $m=3$ and $p=25$, then $\text{k-rank}(\M{X})
\geq 18$ is sufficient for uniqueness. 
The table does not directly depend on $n$; however,
recall that $\MX$ is an $n \times p$ matrix, so $\text{k-rank}(\MX) \leq \min \{n,p\}$.

\begin{table}
  \centering
  \begin{tabular}{rc|*{8}{c}}
\multicolumn{2}{c|}{} & \multicolumn{8}{c}{components ($p$)} \\ 
\multicolumn{2}{c|}{} & 2 & 3 & 4 & 5 & 10 & 25 & 50 & 100 \\ \hline 
\multirow{4}{*}{\rotatebox[origin=c]{90}{order ($m$)}}
& 3 & 2 & 3 & 4 & 4 & 8 & 18 & 34 & 68 \\ 
& 4 & 2 & 3 & 3 & 4 & 6 & 14 & 26 & 51 \\ 
& 5 & 2 & 2 & 3 & 3 & 5 & 11 & 21 & 41 \\ 
& 6 & 2 & 2 & 3 & 3 & 5 & 10 & 18 & 35 \\ 
  \end{tabular}
  \caption{Minimal $\text{k-rank}(\M{X})$ sufficient for uniqueness of symmetric outer product factorization.}
  \label{tab:krank}
\end{table}

The importance of essential uniqueness is that the global solution of the unconstrained problem \Eqn{cp} is the same as for the symmetric problem \Eqn{scp} so long as $\MX$ satisfies \Eqn{unique}.
If we normalize the factors in \Eqn{cp} and, without loss of generality, ignore
the permutation ambiguity, then uniqueness implies, for $k=1,\dots,p$,
\begin{displaymath}
  \Slk = \pm \| \MnC{u}{k}{1} \| \cdots \| \MnC{u}{k}{m} \|
  \qtext{and} 
  \Vxk = \pm \MnC{U}{1}{k}/\| \MnC{u}{k}{1} \| = \cdots = \pm \MnC{U}{m}{k}/\| \MnC{u}{k}{m} \| 
\end{displaymath}
A bit of care must be taken to convert from a solution that ignores
symmetry since it could be the case, e.g., that $\MnC{U}{1}{k} = -
\MnC{U}{2}{k}$. \Alg{symmetrize} gives a simple procedure to
``symmetrize'' a tensor so that the signs align. It also
averages the final sign-aligned factor matrices in case they are not
exactly equal.

The benefit of ignoring symmetry is that we can use existing
software for the CP decomposition. 
The disadvantage is that it requires $m$ times as much storage, 
i.e., it must store the matrices $\Mn{U}{1}$ thru
$\Mn{U}{m}$ rather than just $\MX$.
Moreover, there is no guarantee that the optimization algorithm will find the global minimum.

\begin{algorithm}
  \caption{Symmetrize Kruskal tensor}
  \label{alg:symmetrize}
  Input: CP decomposition defined by 
  $\Mn{U}{1},\dots,\Mn{U}{m}$\\
  Output: Symmetric CP decomposition defined by $\Vl$ and $\MX$
  \begin{algorithmic}[1]
    \For{$k=1,\dots,p$}
    \State $\lambda_k \gets 1$
    \For{$j=1,\dots,m$}
    \State $\eta \gets \| \MnC{U}{j}{k} \|_2$
    \State $\lambda_k \gets \eta \lambda_k$ and $\MnC{U}{j}{k} \gets \MnC{U}{j}{k} / \eta$
    \Comment{Normalize}
    \If{$j> 1$ and $\braket{ \MnC{U}{1}{k}, \MnC{U}{j}{k} } < 0$}
    \State $\lambda_k \gets - \lambda_k$ and $\MnC{U}{j}{k} \gets - \MnC{U}{j}{k} $
    \Comment{Flip $\MnC{U}{j}{k}$ to align with $\MnC{U}{1}{k}$}
    \EndIf
    \EndFor
    \EndFor
    \State $\MX \gets \sum_j \Mn{U}{j}/m$.
  \end{algorithmic}
\end{algorithm}
\section{Optimization formulation for nonnegative symmetric factorization}
\label{sec:nonnegative}

The notion of completely positive tensors has been introduced by Qi,
Xu, and Xu~\cite{QiXuXu13}. It is a natural extension of completely
positive matrices. A nonnegative tensor $\TA \in \STMN$ is
called completely positive if it has a decomposition of the form in
\Eqn{nnscp}.

The formulation is analogous to the unconstrained case, except that there is no $\Vl$ (or equivalently, we constrain $\Vl=1$) and we add nonnegativity constraints.
For given $\TA \in \STMN$, our goal is to find $\MX$ such that \Eqn{nnscp} is satisfied. 
We again assume $p$ is known.
The mathematical program is given by
\begin{displaymath}
  \min f_+(\MX) = \sum_{\Vi \in \mathcal{I}} w_{\Vi} \left( a_{\Vi} - \sum_{k=1}^p (\Vxkm)_{\Vi} \right)^2
  = \sum_{\Vi \in \mathcal{I}} w_{\Vi} \delta_{\Vi}^2 \qtext{s.t.} \M{X} \geq 0.
\end{displaymath}
Choosing $w_{\Vi} = \sigma_{\Vi}$ yields the analogue of $f_1$ whereas $w_{\Vi} = 1$ yields the analogue $f_2$.
The value $\delta_{\Vi}$ is the difference between the model and the tensor at entry $\Vi$.

Using the generic formulation and following \Eqn{grad} without
$\lambda_k$, the gradients are given by
\begin{align*}
  \FD{f_+}{x_{jk}} & = -2 \sum_{\Vc \in \mathcal{C}} 
  c_j \, w_{\Vc} \, \delta_{\Vc} \,
  \left( x_{1k}^{c_1} \cdots x_{jk}^{c_j-1} \cdots x_{nk}^{c_n} \right).
\end{align*}

Our formulation finds the best nonnegative factorization.
Fan and Zhou~\cite{FaZh14} consider the problem of verifying that a
tensor is completely positive.

\section{Numerical results}
\label{sec:results}

For our numerical results, we assume the tensor has underlying
low-rank structure, as is typical in comparisons of numerical methods
for tensor factorization (see, e.g., \cite{ToBr06}). 
Hence, we assume there is some underlying $\Vl^* \in \Real^p$ and
$\MX^* \in \Real^{n \times p}$ to be recovered, where $p$ is lower
than the typical rank. 
The noise-free data tensor is given by
\begin{equation}
  \label{eq:Astar}
  \TA^* = \sum_{k=1}^p \lambda_k^* (\Vxk^* )^m.
\end{equation}
The data tensor $\TA$ may also be contaminated
by noise as controlled by the parameter $\eta \geq 0$, i.e.,
\begin{equation}
  \label{eq:noise}
  \TA = \TA^*  + \eta \frac{\|\TA^*\|}{\|\T{N}\|} \T{N}
  \qtext{where}
  n_{\Vi} \sim \mathcal{N}(0,1).
\end{equation}
Here $\T{N}$ is a noise tensor such that each element is drawn from a
normal distribution, i.e., $n_{\Vi} \sim \mathcal{N}(0,1)$.
The parameters $m$, $n$, $p$ control the size of the problem.  If the
vectors in $\MX^*$ are collinear, then the problem is generally more
difficult \cite{Ki98b,ToBr06}. 

For the $f_1$ objective function in \Eqn{f1}, we calculate the
gradients as specified in \Eqn{f1-grad-alt}. For small problems this
may not be as fast as \Eqn{grad}, but for larger problems it makes a
significant difference in speed, as shown in the results.
For $f_2$, we {precompute} the index set $\mathcal{I}$ as
well as the corresponding monomial representations $\mathcal{C}$ and
multiplicities $\V{\sigma}$. This means that these values need not be
computed each time the objective function and gradient are evaluated. 
The time for this preprocessing is included in the reported runtimes.

All tests were conducted on a laptop with an Intel Dual Core i7-3667U
CPU and 8 GB of RAM, using MATLAB R2013a. For the optimization, unless
otherwise noted, all tests are based on SNOPT, Version 7.2-9
\cite{GiMuSa05,SNOPT7}, using the MATLAB MEX interface. SNOPT default
parameters were used except for the following: Major iteration limit =
10,000, New superbasics limit / Superbasics limit = 999, Major
optimality tolerance = 1e-8.  All tensor computations use the Tensor
Toolbox for MATLAB, Version 2.5 \cite{BaKo06,BaKo07,TTB_Software} as
well as additional codes for symmetric tensors (e.g., to calculate the
index sets) that will be included in the next release.

\subsection{Numerical results on a collection of test problems}
\label{sec:standard}

We consider the impact of the problem formulation resulting from the
choice of objective function and column normalization penalty. The
objective function can \emph{weighted}, based on the standard least
squares formulation denoted by $f_1$ in \Eqn{f1}, or
\emph{unweighted}, which counts each unique entry only once denoted by
$f_2$ in \Eqn{f2}. The column normalization penalty  is either
$\gamma = 0$ (no penalty) or $\gamma = 0.1$. Higher values of $\gamma$
did not change the results.

We test the choices for several test problems
as follows. We consider four sizes:
\begin{itemize}
\item $m=3,n=4,p=2$;
\item  $m=4,n=3,p=5$;
\item  $m=4,n=25,p=3$; and
\item $m=6,n=6,p=4$.
\end{itemize}
In the first case, since $m$ is odd, we have the option to exclude
$\Vl$ from the optimization, but we include it here for
consistency in this set of experiments. 
For each size, we also consider three noise levels: $\eta \in
\set{0,0.01,0.1}$. 

A random instance is created as follows. We generate a \emph{true
solution} defined by $\Vl^* \in \Real^p$ and $\MX^* \in \Real^{n \times
  p}$. The
weight vector has entries selected uniformly from $\{-1,1\}$, i.e.,
\begin{displaymath}
  \Vl^* \in \Real^p \qtext{such that} \lambda_k^* \in \mathcal{U}\{-1,1\}.
\end{displaymath}
The factor matrix is computed by first generating a matrix with random values from the normal distribution, i.e.,
\begin{displaymath}
  \M[\hat]{X}^* \in \Real^{n \times p} 
  \qtext{such that} \hat x_{ik}^* \in \mathcal{N}(0,1),
\end{displaymath}
and then normalizing each column to length one, i.e., $\Vxk^* = \V[\hat]{x}_k^* / \|\V[\hat]{x}_k^*\|=1$.
Finally, given $\Vl^*$ and $\MX^*$, we can compute the tensor
$\T{A}^*$ from \Eqn{Astar} and add noise at the level specified by
$\eta$ per \Eqn{noise}.  For each problem size and  
noise level, we generate ten instances.

For each problem size, we generate five random starting points by choosing
entries of $\MX$ from a Gaussian distribution (no column normalization) and entries of $\Vl$ uniformly at random from $\set{-1,1}$. The same five starting points are used for all problems of that size.

\newcommand{\gammaheader}[1]{\multicolumn{#1}{c|}{$\gamma = 0$} & \multicolumn{#1}{c|}{$\gamma = 0.1$}}
\begin{table}[htbp]
  \centering
\subfloat[\textbf{Relative error}: \textcolor{color1}{runs $\leq
  0.1$}, \textcolor{color2}{instances $\leq  0.1$}, and
\textcolor{color3}{median}.]%
{
  \label{tab:relerr}
  \begin{tabular}{|c@{\,\,\,}c@{\,\,\,}cc|*{4}{>{\color{color1}}r@{/}>{\color{color2}}r@{/}>{\color{color3}}c|}}
    \hline
    \multicolumn{3}{|c}{Size} & Noise & \multicolumn{6}{c|}{Unweighted $f_2$} & \multicolumn{6}{c|}{Weighted $f_1$}\\
    $m$ & $n$ & $p$ & $\eta$ & \gammaheader{3} & \gammaheader{3} \\
    \hline
3 &  4 & 2 & 0.00 & 48 & 10 & 9e-07 & 50 & 10 & 1e-06 & 35 & 10 & 7e-07 & 42 & 10 & 4e-07 \\ 
  &    &   & 0.01 & 43 & 10 & 8e-03 & 46 & 10 & 8e-03 & 32 &  9 & 8e-03 & 39 & 10 & 8e-03 \\ 
  &    &   & 0.10 & 48 & 10 & 8e-02 & 48 & 10 & 8e-02 & 39 & 10 & 8e-02 & 41 & 10 & 8e-02 \\ \hline 
4 &  3 & 5 & 0.00 & 34 &  9 & 5e-02 & 38 & 10 & 6e-03 & 27 & 10 & 9e-02 & 37 & 10 & 4e-02 \\ 
  &    &   & 0.01 & 31 &  9 & 5e-02 & 39 &  9 & 6e-03 & 29 & 10 & 7e-02 & 39 & 10 & 9e-03 \\ 
  &    &   & 0.10 & 36 & 10 & 6e-02 & 39 &  9 & 4e-02 & 38 & 10 & 5e-02 & 40 & 10 & 4e-02 \\ \hline 
4 & 25 & 3 & 0.00 &  6 &  5 & 7e-01 & 40 & 10 & 2e-05 &  4 &  4 & 6e-01 & 16 &  9 & 6e-01 \\ 
  &    &   & 0.01 & 10 &  7 & 7e-01 & 44 & 10 & 1e-02 & 10 &  8 & 6e-01 & 19 & 10 & 6e-01 \\ 
  &    &   & 0.10 & 11 &  7 & 7e-01 & 44 & 10 & 1e-01 & 11 &  6 & 6e-01 & 26 & 10 & 1e-01 \\ \hline 
6 &  6 & 4 & 0.00 & 23 & 10 & 4e-01 & 39 & 10 & 2e-05 &  7 &  5 & 5e-01 & 18 &  9 & 4e-01 \\ 
  &    &   & 0.01 & 15 &  9 & 5e-01 & 40 & 10 & 1e-02 &  9 &  8 & 5e-01 & 25 & 10 & 1e-01 \\ 
  &    &   & 0.10 &  1 &  1 & 5e-01 &  5 &  1 & 1e-01 &  7 &  7 & 5e-01 & 18 & 10 & 3e-01 \\ \hline 
\multicolumn{4}{|c|}{Total}& 306 & 97 & 1e-01 & 472 & 109 & 1e-02 & 248 & 97 & 3e-01 & 360 & 118 & 9e-02 \\ \hline 
  \end{tabular}
}\\
\subfloat[\textbf{Solution score}: \textcolor{color1}{runs $\geq
  0.9$}, \textcolor{color2}{instances $\geq  0.9$}, and
\textcolor{color3}{median}.]%
{
  \label{tab:sc}
  \begin{tabular}{|c@{\,\,\,}c@{\,\,\,}cc|*{4}{>{\color{color1}}r@{/}>{\color{color2}}r@{/}>{\color{color3}}c|}}
\hline
\multicolumn{3}{|c}{Size} & Noise & \multicolumn{6}{c|}{Unweighted $f_2$} & \multicolumn{6}{c|}{Weighted $f_1$}\\
$m$ & $n$ & $p$ & $\eta$ & \gammaheader{3} & \gammaheader{3} \\
\hline
3 &  4 & 2 & 0.00 & 48 & 10 & 1.00 & 50 & 10 & 1.00 & 35 & 10 & 1.00 & 42 & 10 & 1.00 \\ 
  &    &   & 0.01 & 43 & 10 & 0.99 & 46 & 10 & 1.00 & 32 &  9 & 0.99 & 39 & 10 & 1.00 \\ 
  &    &   & 0.10 & 34 &  7 & 0.97 & 34 &  7 & 0.97 & 30 &  8 & 0.93 & 32 &  8 & 0.93 \\ \hline 
4 &  3 & 5 & 0.00 &  3 &  2 & 0.43 &  5 &  3 & 0.64 &  0 &  0 & 0.42 &  6 &  4 & 0.60 \\ 
  &    &   & 0.01 &  0 &  0 & 0.43 &  1 &  1 & 0.54 &  0 &  0 & 0.37 &  2 &  2 & 0.52 \\ 
  &    &   & 0.10 &  0 &  0 & 0.25 &  0 &  0 & 0.45 &  0 &  0 & 0.32 &  1 &  1 & 0.51 \\ \hline 
4 & 25 & 3 & 0.00 &  7 &  6 & 0.55 & 40 & 10 & 1.00 &  9 &  7 & 0.66 & 16 &  9 & 0.67 \\ 
  &    &   & 0.01 & 10 &  7 & 0.53 & 44 & 10 & 1.00 & 11 &  8 & 0.67 & 19 & 10 & 0.67 \\ 
  &    &   & 0.10 & 13 &  7 & 0.51 & 44 & 10 & 1.00 & 16 &  9 & 0.67 & 26 & 10 & 1.00 \\ \hline 
6 &  6 & 4 & 0.00 & 21 & 10 & 0.72 & 38 & 10 & 1.00 &  6 &  4 & 0.72 & 15 &  8 & 0.75 \\ 
  &    &   & 0.01 & 15 &  9 & 0.73 & 40 & 10 & 1.00 &  9 &  8 & 0.67 & 25 & 10 & 0.87 \\ 
  &    &   & 0.10 & 18 &  8 & 0.72 & 32 & 10 & 0.98 &  7 &  7 & 0.73 &
  18 & 10 & 0.74 \\ \hline 
\multicolumn{4}{|c|}{Total}& 212 & 76 & 7e-01 & 374 & 91 & 1e+00 & 155 & 70 & 6e-01 & 241 & 92 & 7e-01 \\ \hline 
  \end{tabular}
}\\
\subfloat[\textbf{Run time}: \textcolor{color1}{mean} and
\textcolor{color2}{standard deviation}.]%
{
  \label{tab:runtime}
  \begin{tabular}{|c@{\,\,\,}c@{\,\,\,}cc|*{4}{>{\color{color1}}r@{$\,\,\pm\,\,$}>{\color{color2}}r|}}
\hline
\multicolumn{3}{|c}{Size} & Noise & \multicolumn{4}{c|}{Unweighted $f_2$} & \multicolumn{4}{c|}{Weighted $f_1$}\\
$m$ & $n$ & $p$ & $\eta$ & \gammaheader{2} & \gammaheader{2} \\
\hline
3 &  4 & 2 & 0.00 & 0.10 & 0.02 & 0.13 & 0.03 & 0.71 & 0.99 & 0.69 & 0.63 \\ 
  &    &   & 0.01 & 0.31 & 0.74 & 0.21 & 0.27 & 0.92 & 1.09 & 1.03 & 1.19 \\ 
  &    &   & 0.10 & 0.12 & 0.08 & 0.17 & 0.28 & 0.92 & 1.47 & 1.02 & 1.48 \\ \hline 
4 &  3 & 5 & 0.00 & 1.19 & 1.81 & 0.75 & 0.87 & 4.67 & 3.29 & 6.28 & 5.44 \\ 
  &    &   & 0.01 & 0.96 & 1.25 & 0.77 & 0.72 & 5.60 & 3.88 & 7.07 & 5.16 \\ 
  &    &   & 0.10 & 1.43 & 1.48 & 0.98 & 0.78 & 6.35 & 4.60 & 6.05 & 4.85 \\ \hline 
4 & 25 & 3 & 0.00 & 41.38 & 20.20 & 54.90 & 13.09 & 3.80 & 2.03 & 5.02 & 1.45 \\ 
  &    &   & 0.01 & 44.96 & 25.17 & 55.69 & 21.31 & 4.42 & 2.34 & 5.47 & 1.99 \\ 
  &    &   & 0.10 & 44.32 & 24.24 & 56.11 & 13.29 & 5.01 & 2.65 & 8.95 & 3.42 \\ \hline 
6 &  6 & 4 & 0.00 & 1.79 & 1.31 & 1.64 & 0.48 & 4.91 & 2.76 & 6.55 & 2.29 \\ 
  &    &   & 0.01 & 1.52 & 0.89 & 1.57 & 0.41 & 7.20 & 5.17 & 8.09 & 2.88 \\ 
  &    &   & 0.10 & 1.57 & 0.74 & 1.76 & 0.73 & 6.16 & 2.88 & 8.56 & 4.08 \\ \hline 
  \end{tabular}
}
  \caption{Results of different formulations for a set of test
    problems. For each size and noise combination, the number of runs
    is fifty and the number of instances is ten (five random starts
    per instance).}
  \label{tab:standard}
\end{table}

For each problem formulation corresponding to a choice for objective function and
for normalization penalty, we do fifty runs, i.e., ten instances with
five random starts each. 
The same instances and starting points are used across all formulations.
The output of each run is a weight vector $\Vl$ and a matrix $\MX$.
\Tab{relerr} compares the relative error which measures the proportion of the
observed data that is explained by the model, i.e., 
\begin{displaymath}
  \text{relative error } = 
  {\left \| \TA - \displaystyle\sum_{k=1}^p \Slk\Pm \Vxk^m  \right\|} /
  { \| \TA \|}.
\end{displaymath}
In the case of no noise, the ideal relative error is zero; otherwise,
we hope for something near the noise level, i.e., $\eta$. In our
experiments, we say a run or instance is \emph{successful} if the relative
error is $\leq$ 0.1.
For each formulation, three values are reported. The first value is the
number of successful runs. Since
we are using five starting points per instance, the second value is
the number of instances such that at least one starting point is successful. 
Finally, the last value is the median
relative error across all fifty runs.
Summary totals are provide in the last line for the 600 runs and 150 instances.
Clearly, $\gamma=0.1$ is superior to $\gamma=0$ in terms of number of
successful runs and instances.
The comparison of unweighted ($f_2$) and weighted ($f_1$) is less
clear cut --- the unweighted formulation is successful for many more
runs overall, but the weighted formulation is successful for more instances
overall. 

\Tab{sc} compares the solution scores which is a measure of
how accurately $\Vl$ and $\MX$ are as compared to $\Vl^*$ and
$\MX^*$. 
Without loss of generality, we assume
both $\MX$ and $\MX^*$ have normalized columns. (If $\| \Vxk \|_2 \neq
1$, then we rescale $\Slk = \Slk \sqrt[m]{ \| \Vxk \|}$ and $\Vxk =
\Vxk / \| \Vxk \|$.)  There is a permutation ambiguity, but we permute
the computed solution so as to maximize the following score:
\begin{displaymath}
  \text{solution score } = \frac{1}{p} \sum_{k=1}^p 
  \left(
    1 - \frac{|\lambda_k - \lambda_k^*|}{ \max\{|\lambda_k|,|\lambda_k^*|\} }
  \right)
  \left| 
    \Vxk^{\Tra} \Vxk^*
  \right|.
\end{displaymath}
A solution score of 1 indicates a perfect match, and we say a run or
instance is \emph{successful} if its solution score is $\geq$ 0.9.
As with the relative error, we report three values. The first value is
the number of runs out of fifty that are successful, the second value is the number
of instances out of ten that are successful (i.e., at least one starting point
was successful), and the third value is the median solution score.
We also report totals for each formulation across the 600 runs and 150 instances.
Consistent with \Tab{relerr}, using $\gamma=0.1$ is more successful
than $\gamma=0$. The unweighted is once again successful for more
runs, but the two methods are nearly tied in terms of number of instances.

Observe in \Tab{sc} that the second size ($m=4,n=3,p=5$) has very low solution
scores despite having good performance in terms of relative
error. This is because the solution may not be unique, i.e., the
k-rank of $\MX^*$ is no more than 3, but the minimum k-rank that is
sufficient for uniqueness is 4 per \Tab{krank}. If the solution is not
unique, then multiple solutions exist and there is no reason to expect
that the particular solution we find will be that one.
For example, a particular instance for $m=4,n=3,p=5$ with $\eta=0$ is defined by
\begin{displaymath}
\V{\lambda}^* = \begin{bmatrix*}[r] 1 \\1 \\1 \\-1 \\-1 \\\end{bmatrix*}
\qtext{and}
\M{X}^* = \begin{bmatrix*}[r] 
-0.3859 & -0.9285 &  0.4922 & -0.1094 &  0.4107 \\
 0.8403 & -0.1678 & -0.6975 &  0.8395 &  0.0308 \\
 0.3807 &  0.3313 & -0.5208 & -0.5322 &  0.9112 \\
\end{bmatrix*}.
\end{displaymath}
The alternate model given by
\begin{displaymath}
\V{\lambda} = \begin{bmatrix*}[r] 1 \\1 \\1 \\-1 \\-1 \\\end{bmatrix*}
\qtext{and}
\M{X} = \begin{bmatrix*}[r] 
-0.7872 &  0.5136 & -0.7809 & -0.1081 &  0.4157 \\
-0.1928 & -0.9150 & -0.0704 &  0.8249 &  0.0387 \\
 0.2039 & -0.5355 &  0.3678 & -0.5477 &  0.9065 \\
\end{bmatrix*}
\end{displaymath}
has a relative error less than $10^{-6}$. The last two columns
generally agree, but the first three do not and the solution score is
only 0.65.
It may be interesting to know that in the matrix case ($m=2$), we
would never compare the computed solution without imposing additional
constraints such as orthogonality.

\Tab{runtime} compares the total runtimes for each method. As with any
nonconvex optimization problem, there is significant variation from
run to run, but we can gain a sense of the general expense for each
method. As a reminder, we computed the gradient in the weighted case
as shown in \Eqn{f1-grad-alt}. If we compute it instead using
\Eqn{grad}, the runtimes for the weighted and unweighted methods are
roughly the same. For size $m=4,n=25,p=3$, the computation in
\Eqn{f1-grad-alt} yields a 5-15X speed improvement because $n$ is
large; otherwise for smaller $n$, the computation in \Eqn{grad} will
generally be faster.

\begin{table}[htbp]
  \centering \begin{tabular}{|c@{\,\,\,}c@{\,\,\,}cc|*{4}{>{\color{color1}}r@{/}>{\color{color2}}r|}}
\hline
\multicolumn{3}{|c}{Size} & Noise & \multicolumn{4}{c|}{Unweighted $f_2$} & \multicolumn{4}{c|}{Weighted $f_1$}\\
$m$ & $n$ & $p$ & $\eta$ & \gammaheader{2} & \gammaheader{2} \\
\hline
3 &  4 & 2 & 0.00 &  0 & 6.73e+00 & 50 & 1.17e-06 &  1 & 3.54e+02 & 50 & 1.87e-06 \\ 
  &    &   & 0.01 &  1 & 3.76e+01 & 50 & 9.11e-05 &  1 & 2.04e+02 & 49 & 9.47e-04 \\ 
  &    &   & 0.10 &  0 & 4.32e+01 & 50 & 9.32e-07 &  0 & 3.16e+02 & 48 & 6.36e-04 \\ \hline 
4 &  3 & 5 & 0.00 &  0 & 7.99e+01 & 46 & 7.34e-03 &  0 & 9.24e+01 & 40 & 9.76e-03 \\ 
  &    &   & 0.01 &  0 & 6.56e+01 & 46 & 2.40e-03 &  0 & 6.47e+01 & 39 & 1.09e-02 \\ 
  &    &   & 0.10 &  0 & 2.76e+02 & 44 & 3.20e-03 &  0 & 1.05e+02 & 39 & 1.71e-02 \\ \hline 
4 & 25 & 3 & 0.00 &  0 & 1.70e+03 & 50 & 2.63e-06 &  0 & 8.36e+02 & 49 & 3.65e-02 \\ 
  &    &   & 0.01 &  0 & 1.76e+03 & 50 & 4.35e-06 &  0 & 6.52e+02 & 50 & 1.17e-06 \\ 
  &    &   & 0.10 &  0 & 1.53e+03 & 50 & 2.99e-06 &  0 & 1.13e+03 & 49 & 6.29e-02 \\ \hline 
6 &  6 & 4 & 0.00 &  1 & 2.44e+01 & 50 & 6.04e-05 &  0 & 6.13e+00 & 50 & 5.46e-05 \\ 
  &    &   & 0.01 &  1 & 2.45e+01 & 50 & 1.74e-05 &  0 & 2.11e+01 & 50 & 3.29e-04 \\ 
  &    &   & 0.10 &  0 & 3.12e+01 & 50 & 3.34e-05 &  0 & 4.42e+01 & 49 & 4.75e-04 \\ \hline 
  \end{tabular}
  \caption{\textbf{Constraint violation}:  \textcolor{color1}{runs  $\leq$ 0.01} and \textcolor{color2}{mean}.}
  \label{tab:cv}
\end{table}

Finally, we briefly consider the impact on $\gamma$ with respect to
the \emph{constraint violation} from \Eqn{gamma}, i.e.,
\begin{displaymath}
  \text{constraint violation } = \sum_{k=1}^p \Parens{ \Vxk^{\Tra}\Vxk - 1 }^2.
\end{displaymath}
In \Tab{cv}, we report the number of runs where the constraint
violation is $\leq$ 0.01 and the mean value.
Recall that the addition of the constraint violation is mainly a
convenience, but it does improve the formulation by eliminating a
manifold of equivalent solutions.

\begin{table}[htbp]
  \centering
\subfloat[\textbf{Relative error}: \textcolor{color1}{runs $\leq
  0.1$}, \textcolor{color2}{instances $\leq  0.1$}, and
\textcolor{color3}{median}.]%
{
  \label{tab:co-relerr}
  \begin{tabular}{|c@{\,\,\,}c@{\,\,\,}cc|*{4}{>{\color{color1}}r@{/}>{\color{color2}}r@{/}>{\color{color3}}c|}}
    \hline
    \multicolumn{3}{|c}{Size} & Noise & \multicolumn{6}{c|}{Unweighted $f_2$} & \multicolumn{6}{c|}{Weighted $f_1$}\\
    $m$ & $n$ & $p$ & $\eta$ & \gammaheader{3} & \gammaheader{3} \\
    \hline
3 &  4 & 2 & 0.00 & 48 & 10 & 9e-07 & 50 & 10 & 1e-06 & 35 & 10 & 7e-07 & 42 & 10 & 4e-07 \\ 
  &    &   & 0.01 & 43 & 10 & 8e-03 & 46 & 10 & 8e-03 & 32 &  9 & 8e-03 & 39 & 10 & 8e-03 \\ 
  &    &   & 0.10 & 48 & 10 & 8e-02 & 48 & 10 & 8e-02 & 39 & 10 & 8e-02 & 41 & 10 & 8e-02 \\ \hline 
4 & 25 & 3 & 0.00 &  8 &  7 & 7e-01 & 40 & 10 & 2e-05 &  4 &  4 & 6e-01 & 16 &  9 & 6e-01 \\ 
  &    &   & 0.01 &  9 &  6 & 7e-01 & 44 & 10 & 1e-02 & 10 &  8 & 6e-01 & 19 & 10 & 6e-01 \\ 
  &    &   & 0.10 & 11 &  6 & 7e-01 & 44 & 10 & 1e-01 & 11 &  6 & 6e-01 & 26 & 10 & 1e-01 \\ \hline 
6 &  6 & 4 & 0.00 & 11 &  6 & 3e-01 & 26 &  9 & 2e-02 &  2 &  2 & 3e-01 &  6 &  5 & 2e-01 \\ 
  &    &   & 0.01 & 17 & 10 & 2e-01 & 30 &  9 & 1e-02 &  4 &  3 & 3e-01 & 16 & 10 & 2e-01 \\ 
  &    &   & 0.10 &  3 &  2 & 2e-01 &  8 &  2 & 1e-01 & 11 &  8 & 2e-01 & 19 & 10 & 1e-01 \\ \hline 
  \end{tabular}
}\\
\subfloat[\textbf{Solution score}: \textcolor{color1}{runs $\geq
  0.9$}, \textcolor{color2}{instances $\geq  0.9$}, and
\textcolor{color3}{median}.]%
{
  \label{tab:co-sc}
  \begin{tabular}{|c@{\,\,\,}c@{\,\,\,}cc|*{4}{>{\color{color1}}r@{/}>{\color{color2}}r@{/}>{\color{color3}}c|}}
\hline
\multicolumn{3}{|c}{Size} & Noise & \multicolumn{6}{c|}{Unweighted $f_2$} & \multicolumn{6}{c|}{Weighted $f_1$}\\
$m$ & $n$ & $p$ & $\eta$ & \gammaheader{3} & \gammaheader{3} \\
\hline
3 &  4 & 2 & 0.00 & 48 & 10 & 9e-07 & 50 & 10 & 1e-06 & 35 & 10 & 7e-07 & 42 & 10 & 4e-07 \\ 
  &    &   & 0.01 & 43 & 10 & 8e-03 & 46 & 10 & 8e-03 & 32 &  9 & 8e-03 & 39 & 10 & 8e-03 \\ 
  &    &   & 0.10 & 48 & 10 & 8e-02 & 48 & 10 & 8e-02 & 39 & 10 & 8e-02 & 41 & 10 & 8e-02 \\ \hline 
4 & 25 & 3 & 0.00 &  8 &  7 & 7e-01 & 40 & 10 & 2e-05 &  4 &  4 & 6e-01 & 16 &  9 & 6e-01 \\ 
  &    &   & 0.01 &  9 &  6 & 7e-01 & 44 & 10 & 1e-02 & 10 &  8 & 6e-01 & 19 & 10 & 6e-01 \\ 
  &    &   & 0.10 & 11 &  6 & 7e-01 & 44 & 10 & 1e-01 & 11 &  6 & 6e-01 & 26 & 10 & 1e-01 \\ \hline 
6 &  6 & 4 & 0.00 & 11 &  6 & 3e-01 & 26 &  9 & 2e-02 &  2 &  2 & 3e-01 &  6 &  5 & 2e-01 \\ 
  &    &   & 0.01 & 17 & 10 & 2e-01 & 30 &  9 & 1e-02 &  4 &  3 & 3e-01 & 16 & 10 & 2e-01 \\ 
  &    &   & 0.10 &  3 &  2 & 2e-01 &  8 &  2 & 1e-01 & 11 &  8 & 2e-01 & 19 & 10 & 1e-01 \\ \hline 
  \end{tabular}
}%
  \caption{Results of different formulations for ``collinear test
    problems. For each size and noise combination, the number of runs
    is fifty and the number of instances is ten (five random starts
    per instance).}
  \label{tab:co}
\end{table}

\Tab{co} shows results for more difficult test problems where
$\MX^*$ has collinear normalized columns,i.e.,
$(\Vxk^*)^{\Tra}\V{x}_{\ell}^*= 0.9$ for all $k \neq \ell$ with
$k,\ell \in \{1,\dots,p\}$. The procedure for generating the collinear
columns is described by Tomasi and Bro \cite{ToBr06}.
The setup is the same as in the previous subsection except for the
change in how we generate $\MX^*$ and the omission of size
$m=4,n=3,p=5$ (since the procedure we are using does not allow $p>n$).
The results in \Tab{co} are are analogous to those in
\Tab{standard}. We omit the runtimes since they are similar. Although
fewer runs are successful, the number of instances solved is similar.

From these results, we have a sense that the symmetric factorization
problem can be solved using standard optimization techniques.
Because the problems are nonconvex, multiple starting points are
needed to improve
the odds of finding a global minimizer.  Our results also indicate
that it is helpful
to add a penalty to remove the scaling ambiguity; otherwise, with no
penalty, the Jacobian at the solution is singular which seems to have
a negative impact on the solution quality. 
\subsection{Ignoring symmetry}

As noted previously, Carroll and Chang~\cite{CaCh70} ignored symmetry
with the idea that it may not be required.  Ideally, the solution that
is computed by a standard method, like CP-ALS \cite{FaBrHo03,KoBa09}
or CP-OPT \cite{AcDuKo11}, will be symmetric up to scaling.

Using the same problems from \Tab{standard}, we apply CP-ALS (as
implemented in the Tensor Toolbox), followed by \Alg{symmetrize} to
symmetrize the solution.  Three of the four sizes generically satisfy
the sufficient uniqueness condition in \Eqn{unique}.
\begin{itemize}
\item For $m=3$ and $p=2$, we require $\text{k-rank}(\MX^*) \geq
  2$. Since $\MX^*$ is an $n \times p$ matrix with $n=4$ whose columns
  are randomly generated, $\text{k-rank}(\MX^*) = 2$ with probability
  1.
\item For $m=4$ and $p=5$, we require $\text{k-rank}(\MX^*) \geq
  4$. Since $\MX^*$ is an $n \times p$ matrix with $n=3$, it \emph{cannot}
  satisfy the condition because $\text{k-rank}(\MX^*) \leq
  \text{rank}(\MX^*) \leq \min\set{n,p}=3$. Hence, the solutions may not be unique, and
  an example of a non-unique solution is provided in the previous
  subsection. 
\item For $m=4$ and $p=3$, we require $\text{k-rank}(\MX^*) \geq
  3$. Since $\MX^*$ is an $n \times p$ matrix with $n=25$ whose columns
  are randomly generated, $\text{k-rank}(\MX^*) = 3$ with probability
  1.
\item For $m=6$ and $p=6$, we require $\text{k-rank}(\MX^*) \geq
  3$. Since $\MX^*$ is an $n \times p$ matrix with $n=4$ whose columns
  are randomly generated, $\text{k-rank}(\MX^*) = 4$ with probability
  1.
\end{itemize}
\Tab{unsym} shows the results, which are analogous to those in
\Tab{standard}.
CP-ALS with symmetrization is highly competitive. In terms of the
relative error, its total number of 442 successful runs is near the high
of 472 for the symmetric optimization methods; likewise, it has 116
successful instances versus 118 for symmetric optimization. Its scores
are not as impressive in terms of the solution score, though this is mainly a
problem for the size $m=4,n=3,p=5$, as expected due to lack of symmetry. The major advantage
of CP-ALS is runtime, where it is typically ten times faster or more.
Despite the fact that CP-ALS may not find a symmetric solution, using a
standard CP solution procedure followed by symmetrization is indeed an
effective approach in many situations.

\begin{table}[htbp]
  \centering
  \begin{tabular}
{|c@{\,\,\,}c@{\,\,\,}c|c|
>{\color{color1}}r@{/}
>{\color{color2}}r@{/}
>{\color{color3}}r@{/}
>{\color{color4}}r|
>{\color{color1}}r@{/}
>{\color{color2}}r@{/}
>{\color{color3}}r|
>{\color{color1}}r@{$\,\,\pm\,\,$}
>{\color{color2}}r|}
    \hline
    \multicolumn{3}{|c}{Size} & Noise & \multicolumn{9}{c|}{CP-ALS + Symmetrization}\\
    $m$ & $n$ & $p$ & $\eta$ 
& \multicolumn{4}{c|}{Relative Error}
& \multicolumn{3}{c|}{Soln.\@ Score}
& \multicolumn{2}{c|}{Runtime}
\\
    \hline
3 &  4 & 2 & 0.00 & 44 & 10 & 2e-04 & 2e-04 & 44 & 10 & 1.00 & 0.07 & 0.05 \\ 
  &    &   & 0.01 & 42 & 10 & 8e-03 & 8e-03 & 40 & 10 & 0.99 & 0.06 & 0.04 \\ 
  &    &   & 0.10 & 47 & 10 & 8e-02 & 8e-02 & 40 &  9 & 0.97 & 0.04 & 0.04 \\ \hline 
4 &  3 & 5 & 0.00 & 39 & 10 & 3e-02 & 3e-02 &  3 &  3 & 0.63 & 0.23 & 0.08 \\ 
  &    &   & 0.01 & 36 & 10 & 4e-02 & 2e-02 &  1 &  1 & 0.57 & 0.22 & 0.10 \\ 
  &    &   & 0.10 & 37 & 10 & 4e-02 & 4e-02 &  0 &  0 & 0.59 & 0.21 & 0.09 \\ \hline 
4 & 25 & 3 & 0.00 & 37 &  9 & 9e-07 & 1e-06 & 37 &  9 & 1.00 & 0.07 & 0.04 \\ 
  &    &   & 0.01 & 44 & 10 & 1e-02 & 1e-02 & 44 & 10 & 1.00 & 0.07 & 0.03 \\ 
  &    &   & 0.10 & 46 & 10 & 1e-01 & 1e-01 & 46 & 10 & 1.00 & 0.07 & 0.03 \\ \hline 
6 &  6 & 4 & 0.00 & 29 &  9 & 3e-04 & 3e-04 & 26 &  8 & 1.00 & 0.13 & 0.11 \\ 
  &    &   & 0.01 & 18 &  8 & 5e-01 & 5e-01 & 18 &  8 & 0.73 & 0.08 & 0.06 \\ 
  &    &   & 0.10 & 23 & 10 & 4e-01 & 4e-01 & 23 & 10 & 0.74 & 0.09 & 0.07 \\ \hline 
\multicolumn{4}{|c|}{Total} 
                  & 442 & 116 & 5e-02 & 4e-02 & 322 & 88 & 1e+00 & 
\multicolumn{2}{c|}{}\\ \hline 
  \end{tabular}
  \caption{Results of CP-ALS plus symmetrization on test problems from \Tab{standard}.
    \textbf{Relative error}: 
    \textcolor{color1}{runs $\leq 0.1$}, 
    \textcolor{color2}{instances $\leq  0.1$}, 
    \textcolor{color3}{median symmetrized}, and
    \textcolor{color4}{median unsymmetrized}.
    \textbf{Solution score}: 
    \textcolor{color1}{runs $\geq  0.9$}, 
    \textcolor{color2}{instances $\geq  0.9$}, and
    \textcolor{color3}{median}    
    \textbf{Runtime}: \textcolor{color1}{mean} and
    \textcolor{color2}{standard deviation}.    
  }
  \label{tab:unsym}
\end{table}

\subsection{Sparsity penalty for rank determination}

In Example 5.5(i) of \cite{Ni14}, Nie considers an method for
determining the rank of a tensor. The example tensor is of order $m=4$
and defined by
\begin{displaymath}
  \Vl^* =
  \begin{bmatrix}
    676 \\ 196
  \end{bmatrix}
  \qtext{and}
  \MX^* = 
  \begin{bmatrix*}[r]
    0 & 3/\sqrt{14} \\
    1 / \sqrt{26} & 2 / \sqrt{14} \\
    -5 / \sqrt{26} & -1/\sqrt{14}
  \end{bmatrix*} = 
  \begin{bmatrix*}[r]
   0.00 &    0.80 \\
   0.20 &    0.53 \\
  -0.98 &   -0.27 \\
  \end{bmatrix*}.
\end{displaymath}
Using our optimization approach with $f_1$ and $\gamma=0.1$, we impose
the approximate $\ell_1$ penalty of the form suggested by
\cite{ScFuRo07}, using $\alpha=10$ and $\beta=0.1$ to arrive at the
following result:
\begin{displaymath}
  \Vl = 
  \begin{bmatrix*}[r]
    675.998\\ 195.965\\ 0.001\\ 0.001\\ 0.001\\ 0.001\\ 
  \end{bmatrix*}
  \qtext{and}
  \MX = 
  \begin{bmatrix*}[r]
  -0.00 &    0.80 &   -0.80 &    0.80 &   -0.79 &   -0.02 \\
  -0.20 &    0.53 &   -0.53 &    0.54 &   -0.55 &   -0.26 \\
   0.98 &   -0.27 &    0.27 &   -0.25 &    0.27 &    0.97 \\
  \end{bmatrix*}.
\end{displaymath}
We calculate the similarity score as described previously, selecting
the two components that yield the best match for a score of 0.999865. The
calculation takes approximately 2 seconds.
Using $\alpha=1000$ causes numerical blow-up, but $\alpha=100$ or
$\alpha=1$ work nearly as well as $\alpha=10$, i.e., the solution
score is 0.9998 (with $\beta=0.1$). Likewise, varying $\beta$ has
little impact on the solution quality (with $\alpha=10$).

Using the same penalty parameters ($\alpha=10$ and $\beta=0.1$), we construct 10 instances of
problems of size $m=4$, $n=3$, and $p=2$ for each noise level $\eta \in
\set{0,0.01,0.1}$. We use a solution with 3 components but once again
apply the sparsity penalty, using the same parameters as above. We use
five random starts per instance. The results as shown in
\Tab{l1penalty}. The second column shows the number of instances (out
of 10) where the solution score was $\geq 0.9$, and the third column
is the total number of runs that are successful (out of 50) for which this
condition was satisfied. The fourth column shows that median relative
error, and the last column shows the mean and standard deviation of
the runtime. 
In the noise-free case, the correct solution is found in every run.
For $\eta = 0.01$, the correct solution is obtained for 9 out of 10 instances. 
For $\eta=0.1$, the problem is only solved to the desired accuracy in
4 out of 10 instances.

\begin{table}
  \centering
  \begin{tabular}{|*{5}{c|}}
\hline
& \multicolumn{2}{c|}{Soln. Score $\geq$ 0.9} & Median & Runtime \\
Noise & Instances & Runs & Rel.\@ Error & (mean $\pm$ std.)\\
\hline
$\eta = 0.00$& 10 & 50& 3.4e-04& $0.60 \pm 0.20$\\ \hline 
$\eta = 0.01$& 9 & 45& 6.9e-03& $0.54 \pm 0.19$\\ \hline 
$\eta = 0.10$& 4 & 19& 6.4e-02& $0.47 \pm 0.13$\\ \hline 
  \end{tabular}
  \caption{Impact of sparsity penalty for problems of size $m=4$, $n=3$, and $p=2$ with a solution that has $p=3$.}
  \label{tab:l1penalty}
\end{table}

Alas, the penalty approach is a heuristic; forthcoming work
\cite{woody} will use statistical validation to select the rank.

\subsection{Nonnegative factorization}

Lastly, we consider the problem of nonnegative factorization. 
We use the same problem setup as in \Sec{standard} with the exception
that we set all entries $\Vl^*$ equal to one  and choose entries of
$\MX^*$ to be uniform on $[0,1]$, i.e., $x_{ij}^* \in \mathcal{U}[0,1]$.
The optimization formulation excludes $\Vl$, so there is no penalty on
the columns norms of $\MX$ ($\gamma = 0$). We add bound constraints
that all entries of $\MX$ are nonnegative.
We compare only the weighted and unweighted formulations.
\Tab{nonneg} shows the results, which are analogous to
\Tab{standard}. There is little difference between the two
formulations, except runtimes as discussed previously.

\begin{table}[htbp]
  \centering
  \begin{tabular}{|c@{\,\,}c@{\,\,}c@{\,\,}c|
      *{2}{>{\color{color1}}r@{/}>{\color{color2}}r@{/}>{\color{color3}}c|}
      *{2}{>{\color{color1}}r@{/}>{\color{color2}}r@{/}>{\color{color3}}r|}
      *{2}{r|}
    }
    \hline
    \multicolumn{3}{|c}{Size} & Noise 
    & \multicolumn{6}{c|}{Relative Error} 
    & \multicolumn{6}{c|}{Solution Score} 
    & \multicolumn{2}{c|}{Runtime} \\
    $m$ & $n$ & $p$ & $\eta$ 
& \multicolumn{3}{c|}{Unweighted} & \multicolumn{3}{c|}{Weighted} 
& \multicolumn{3}{c|}{Unw.} & \multicolumn{3}{c|}{Wei.} 
& \multicolumn{1}{c|}{Unw.} & \multicolumn{1}{c|}{Wei.} \\
    \hline
3 &  4 & 2 & 0.00 & 50 & 10 & 4e-07 & 49 & 10 & 1e-07 & 48 & 10 & 1.00 & 49 & 10 & 1.00 & 0.08 & 0.48 \\ 
  &    &   & 0.01 & 50 & 10 & 8e-03 & 50 & 10 & 8e-03 & 30 &  6 & 0.95 & 29 &  6 & 0.93 & 0.11 & 0.53 \\ 
  &    &   & 0.10 & 50 & 10 & 8e-02 & 50 & 10 & 7e-02 & 10 &  2 & 0.71 & 10 &  2 & 0.77 & 0.07 & 0.31 \\ \hline 
4 &  3 & 5 & 0.00 & 50 & 10 & 1e-04 & 50 & 10 & 3e-05 &  2 &  1 & 0.72 &  8 &  6 & 0.76 & 0.30 & 4.64 \\ 
  &    &   & 0.01 & 50 & 10 & 3e-03 & 50 & 10 & 3e-03 &  0 &  0 & 0.56 &  0 &  0 & 0.60 & 0.40 & 4.72 \\ 
  &    &   & 0.10 & 50 & 10 & 5e-02 & 50 & 10 & 4e-02 &  0 &  0 & 0.58 &  0 &  0 & 0.56 & 0.21 & 2.53 \\ \hline 
4 & 25 & 3 & 0.00 & 50 & 10 & 1e-08 & 50 & 10 & 1e-08 & 50 & 10 & 1.00 & 50 & 10 & 1.00 & 20.28 & 2.30 \\ 
  &    &   & 0.01 & 50 & 10 & 1e-02 & 50 & 10 & 1e-02 & 50 & 10 & 1.00 & 50 & 10 & 1.00 & 21.59 & 2.38 \\ 
  &    &   & 0.10 & 50 & 10 & 1e-01 & 50 & 10 & 1e-01 & 50 & 10 & 1.00 & 50 & 10 & 1.00 & 21.58 & 2.27 \\ \hline 
6 &  6 & 4 & 0.00 & 50 & 10 & 2e-07 & 50 & 10 & 1e-08 & 50 & 10 & 1.00 & 50 & 10 & 1.00 & 0.91 & 3.17 \\ 
  &    &   & 0.01 & 50 & 10 & 1e-02 & 50 & 10 & 1e-02 & 50 & 10 & 0.99 & 50 & 10 & 0.99 & 0.88 & 2.70 \\ 
  &    &   & 0.10 &  8 &  2 & 1e-01 & 47 & 10 & 1e-01 & 18 &  4 & 0.77 & 15 &  3 & 0.84 & 0.86 & 2.64 \\ \hline 
  \end{tabular}
  \caption{Results of nonnegative optimization on test problems.
    \textbf{Relative error}: 
    \textcolor{color1}{runs $\leq 0.1$}, 
    \textcolor{color2}{instances $\leq  0.1$}, 
    \textcolor{color3}{median}.
    \textbf{Solution score}: 
    \textcolor{color1}{runs $\geq  0.9$}, 
    \textcolor{color2}{instances $\geq  0.9$}, and
    \textcolor{color3}{median}.    
    \textbf{Runtime}: mean.    
  }
  \label{tab:nonneg}
\end{table}

\section{Conclusions and future challenges}
\label{sec:conclusions}

We consider straightforward optimization formulations for
real-valued symmetric and nonnegative symmetric tensor decompositions. 
These methods can be used as a baselines for comparison as new methods
are developed. 
In particular, these methods should be useful for larger problems with
inherent low-rank structure. For instance, the
size $m=4$ and $n=25$ is larger in terms of dimension than most other
symmetric tensor decomposition 
problems in the literature, though other works consider larger values
of $p$ \cite{Ni14}. 
Furthermore, we consider noise-contaminated problems, which may be
problematic for algebraic methods.

Although the symmetric and nonnegative symmetric tensor decomposition 
problems are nonconvex, these numerical optimization
approaches are effective at recovering the known solution in our
experiments, especially when we use multiple random starting points. 
These optimization formulations can be adapted to the case of partial
symmetries. 
Moreover, we show that if the solution is essentially unique (and the
optimization method finds a global minima), then symmetry need not be
directly enforced by the optimization method. In this case, efficient
tools for the nonsymmetric CP problem may be employed directly.

We expect many further improvements, including different optimization
formulations that exploit structure and consideration of other
optimization methods.

\begin{acknowledgements}
  The anonymous referees provided extremely useful feedback that has
  greatly improved the manuscript.
This material is based upon work supported by the U.S. Department of
Energy, Office of Science, Office of Advanced Scientific Computing
Research, Applied Mathematics program.
Sandia National Laboratories is a multi-program laboratory managed and
operated by Sandia Corporation, a wholly owned subsidiary of Lockheed
Martin Corporation, for the U.S. Department of Energy's National
Nuclear Security Administration under contract DE--AC04--94AL85000. 
\end{acknowledgements}

\bibliographystyle{spmpsci}


\end{document}